\documentstyle[12pt,epsf]{amlt}
\begin{document}
Content-Length: 21789
X-Lines: 501
Status: RO

% Authors, please ignore this section.
%-------------- Publisher's entries --------------------

%\annalsline{0}{1991}
%\received{??}
%\revised{??}
\startingpage{1}

% Authors: Please start here:
%--------------- Author macros ---------------

%\input{amsmac}
%(Optional)
% Please enter all author-written macros 
% used in the body of the paper here:
% i.e.,
% \def\CC{\rm C\!\!\!\!I\,}, etc.

\def\descriptionlabel#1{\normalshape\rm\kern-.5em#1}

\newcommand{\mycirc}{\hbox{\raise.25ex\hbox{${\scriptstyle
\circ}$}}\ }
\newcommand{\ig}{\ell}
\newcommand{\ttot}{\tau_{\mbox{\scriptsize total}}}
\renewcommand{\varphi}{\wp}
\newcommand{\R}{\Bbb{R}} %added by cdh
\newcommand{\h}{\bold{H}} %added by cdh; latex already
%uses \H for something
\newcommand{\s}{\bold{S}} %added by cdh; latex already
%uses \S for something
\newcommand{\E}{\bold{E}} %added by cdh
\newcommand{\p}{{\cal P}}
\newcommand{\m}{{\cal M}}
\newcommand{\thefig}{\thefigure}
\newcommand{\bd}[1]{{\bf #1}}  % \bd{stuff} puts stuff in
%bold
%%\newcommand{\pf}[1]{{\smallskip\noindent\bf #1}}  % for
%proof headings
\newcommand{\fn}[1]{{\tiny #1}}
\newcommand{\cp}[1]{{\sc #1}}  % \cp{stuff} puts stuff in
%small caps
\newcommand{\ic}[1]{{\it #1}}  % \ic{stuff} puts stuff in
%italics
\newcommand{\es}[1]{{\em #1}}  % \es{stuff} puts stuff in
			       % italics/roman toggle
\newcommand{\ro}[1]{\mathrm{#1}}  % \ro{stuff} puts stuff
%in roman

\newcommand{\Kext}{K_{\mbox{\scriptsize ext}}}
\def\QED{\enddemo}
\newcommand{\sphi}{$\bigcirc_\infty$}  % sphere at
%hyperbolic infinity BROKEN!!
\newcommand{\ij}{\ro{InjRad}}  % injectivity radius
\let\mcol=\multicolumn
\def\0{\kern.5em}		% width of a digit
\newcommand{\po}{{\cal P}}  %polar map (script P)
\newcommand{\poi}{{\cal P}^{- 1}}  %inverse polar map
%(script P^-1)
%%\newcommand{\implies}{\Rightarrow}
\def\di{\mathop {\rm dist}\nolimits}     %distance
\def\det{\mathop {\rm det}\nolimits}     %det
\def\tr{\mathop {\rm tr}\nolimits}     %trace
\def\sys{\mathop {\rm sys}\nolimits}     
\def\length{\mathop {\rm length}\nolimits}
\def\Area{\mathop {\rm Area}\nolimits}   
\def\grad{\mathop {\rm grad}\nolimits}     %gradient
\newcommand{\len}{\ell}         %length
\def\arctanh{\mathop {\rm arctanh}\nolimits}
\def\arcsec{\mathop {\rm arcsec}\nolimits}
\def\sign{\mathop {\rm sign}\nolimits}
\def\Aut{\mathop {\rm Aut}\nolimits}
\def\Isom{\mathop {\mbox{\rm Isom}}\nolimits}
\newcommand{\degree}{^{\circ}}
\newcommand{\degrees}{^{\circ}}
\newcommand{\arccosh}{\ro{arccosh}}
\def\arccosh{\mathop {\rm arccosh}\nolimits}
\def\arcsinh{\mathop {\rm arcsinh}\nolimits}
\def\diam{\mathop {\rm diam}\nolimits}
\def\id{\mathop {\rm id}\nolimits}
\newcommand{\Bd}{{\partial }} %added by cdh
\newcommand{\ip}[3]{\langle #1 \rangle_{#2 , #3}}   %inner
%product (new)
\newcommand{\nr}[3]{\| #1 \|_{#2 , #3}}  %norm (new)
\newcommand{\FIG}[1]{\vspace{ #1 }}   % puts #1 (which is
%``1.5in'' or
				      % something) blank space for
				      % pastein

\newcommand{\leb}[1]{} %for uncited labels
\newcommand{\lab}[1]{\label{#1}}
\newcommand{\rf}[1]{\ref{#1}}
\newcommand{\numb}[1]{(\ref{#1})}
\newcommand{\thm}[2]{{\proclaim{#1} #2 \endproclaim}}
\newcommand{\ie}[1]{{\index{#1}}}
\newcommand{\mpf}[1]{\marginpar{\fn#1}}

% This is symbols.tex
% the symbols not available in plain TeX are constructed
% by overprinting some characters

\def\sun{{\hbox{$\odot$}}}
\def\la{\mathrel{\mathchoice {\vcenter{\offinterlineskip\halign{\hfil
$\displaystyle##$\hfil\cr<\cr\noalign{\vskip1.5pt}\sim\cr}}}
{\vcenter{\offinterlineskip\halign{\hfil$\textstyle##$\hfil\cr<\cr
\noalign{\vskip1.0pt}\sim\cr}}}
{\vcenter{\offinterlineskip\halign{\hfil$\scriptstyle##$\hfil\cr<\cr
\noalign{\vskip0.5pt}\sim\cr}}}
{\vcenter{\offinterlineskip\halign{\hfil$\scriptscriptstyle##$\hfil
\cr<\cr\noalign{\vskip0.5pt}\sim\cr}}}}}
\def\ga{\mathrel{\mathchoice {\vcenter{\offinterlineskip\halign{\hfil
$\displaystyle##$\hfil\cr>\cr\noalign{\vskip1.5pt}\sim\cr}}}
{\vcenter{\offinterlineskip\halign{\hfil$\textstyle##$\hfil\cr>\cr
\noalign{\vskip1.0pt}\sim\cr}}}
{\vcenter{\offinterlineskip\halign{\hfil$\scriptstyle##$\hfil\cr>\cr
\noalign{\vskip0.5pt}\sim\cr}}}
{\vcenter{\offinterlineskip\halign{\hfil$\scriptscriptstyle##$\hfil
\cr>\cr\noalign{\vskip0.5pt}\sim\cr}}}}}
\def\sq{\hbox{\rlap{$\sqcap$}$\sqcup$}}
\def\degr{\hbox{$^\circ$}}
\def\arcmin{\hbox{$^\prime$}}
\def\arcsec{\hbox{$^{\prime\prime}$}}
\def\utw{\smash{\rlap{\lower5pt\hbox{$\sim$}}}}
\def\udtw{\smash{\rlap{\lower6pt\hbox{$\approx$}}}}
\def\fd{\hbox{$.\!\!^{\rm d}$}}
\def\fh{\hbox{$.\!\!^{\rm h}$}}
\def\fm{\hbox{$.\!\!^{\rm m}$}}
\def\fs{\hbox{$.\!\!^{\rm s}$}}
\def\fdg{\hbox{$.\!\!^\circ$}}
\def\farcm{\hbox{$.\mkern-4mu^\prime$}}
\def\farcs{\hbox{$.\!\!^{\prime\prime}$}}
\def\fp{\hbox{$.\!\!^{\scriptscriptstyle\rm p}$}}
\def\getsto{\mathrel{\mathchoice {\vcenter{\offinterlineskip
\halign{\hfil$\displaystyle##$\hfil\cr\gets\cr\to\cr}}}
{\vcenter{\offinterlineskip\halign{\hfil$\textstyle##$\hfil\cr
\gets\cr\to\cr}}}
{\vcenter{\offinterlineskip\halign{\hfil$\scriptstyle##$\hfil\cr
\gets\cr\to\cr}}}
{\vcenter{\offinterlineskip\halign{\hfil$\scriptscriptstyle##$\hfil\cr
\gets\cr\to\cr}}}}}
\def\cor{\mathrel{\mathchoice {\hbox{$\widehat=$}}{\hbox{$\widehat=$}}
{\hbox{$\scriptstyle\hat=$}}
{\hbox{$\scriptscriptstyle\hat=$}}}}
\def\grole{\mathrel{\mathchoice {\vcenter{\offinterlineskip\halign{\hfil
$\displaystyle##$\hfil\cr>\cr\noalign{\vskip-1.5pt}<\cr}}}
{\vcenter{\offinterlineskip\halign{\hfil$\textstyle##$\hfil\cr
>\cr\noalign{\vskip-1.5pt}<\cr}}}
{\vcenter{\offinterlineskip\halign{\hfil$\scriptstyle##$\hfil\cr
>\cr\noalign{\vskip-1pt}<\cr}}}
{\vcenter{\offinterlineskip\halign{\hfil$\scriptscriptstyle##$\hfil\cr
>\cr\noalign{\vskip-0.5pt}<\cr}}}}}
\def\lid{\mathrel{\mathchoice {\vcenter{\offinterlineskip\halign{\hfil
$\displaystyle##$\hfil\cr<\cr\noalign{\vskip1.5pt}=\cr}}}
{\vcenter{\offinterlineskip\halign{\hfil$\textstyle##$\hfil\cr<\cr
\noalign{\vskip1pt}=\cr}}}
{\vcenter{\offinterlineskip\halign{\hfil$\scriptstyle##$\hfil\cr<\cr
\noalign{\vskip0.5pt}=\cr}}}
{\vcenter{\offinterlineskip\halign{\hfil$\scriptscriptstyle##$\hfil\cr
<\cr\noalign{\vskip0.5pt}=\cr}}}}}
\def\gid{\mathrel{\mathchoice {\vcenter{\offinterlineskip\halign{\hfil
$\displaystyle##$\hfil\cr>\cr\noalign{\vskip1.5pt}=\cr}}}
{\vcenter{\offinterlineskip\halign{\hfil$\textstyle##$\hfil\cr>\cr
\noalign{\vskip1pt}=\cr}}}
{\vcenter{\offinterlineskip\halign{\hfil$\scriptstyle##$\hfil\cr>\cr
\noalign{\vskip0.5pt}=\cr}}}
{\vcenter{\offinterlineskip\halign{\hfil$\scriptscriptstyle##$\hfil\cr
>\cr\noalign{\vskip0.5pt}=\cr}}}}}
\def\sol{\mathrel{\mathchoice {\vcenter{\offinterlineskip\halign{\hfil
$\displaystyle##$\hfil\cr\sim\cr\noalign{\vskip-0.2mm}<\cr}}}
{\vcenter{\offinterlineskip
\halign{\hfil$\textstyle##$\hfil\cr\sim\cr<\cr}}}
{\vcenter{\offinterlineskip
\halign{\hfil$\scriptstyle##$\hfil\cr\sim\cr<\cr}}}
{\vcenter{\offinterlineskip
\halign{\hfil$\scriptscriptstyle##$\hfil\cr\sim\cr<\cr}}}}}
\def\sog{\mathrel{\mathchoice {\vcenter{\offinterlineskip\halign{\hfil
$\displaystyle##$\hfil\cr\sim\cr\noalign{\vskip-0.2mm}>\cr}}}
{\vcenter{\offinterlineskip
\halign{\hfil$\textstyle##$\hfil\cr\sim\cr>\cr}}}
{\vcenter{\offinterlineskip
\halign{\hfil$\scriptstyle##$\hfil\cr\sim\cr>\cr}}}
{\vcenter{\offinterlineskip
\halign{\hfil$\scriptscriptstyle##$\hfil\cr\sim\cr>\cr}}}}}
\def\lse{\mathrel{\mathchoice {\vcenter{\offinterlineskip\halign{\hfil
$\displaystyle##$\hfil\cr<\cr\noalign{\vskip1.5pt}\simeq\cr}}}
{\vcenter{\offinterlineskip\halign{\hfil$\textstyle##$\hfil\cr<\cr
\noalign{\vskip1pt}\simeq\cr}}}
{\vcenter{\offinterlineskip\halign{\hfil$\scriptstyle##$\hfil\cr<\cr
\noalign{\vskip0.5pt}\simeq\cr}}}
{\vcenter{\offinterlineskip
\halign{\hfil$\scriptscriptstyle##$\hfil\cr<\cr
\noalign{\vskip0.5pt}\simeq\cr}}}}}
\def\gse{\mathrel{\mathchoice {\vcenter{\offinterlineskip\halign{\hfil
$\displaystyle##$\hfil\cr>\cr\noalign{\vskip1.5pt}\simeq\cr}}}
{\vcenter{\offinterlineskip\halign{\hfil$\textstyle##$\hfil\cr>\cr
\noalign{\vskip1.0pt}\simeq\cr}}}
{\vcenter{\offinterlineskip\halign{\hfil$\scriptstyle##$\hfil\cr>\cr
\noalign{\vskip0.5pt}\simeq\cr}}}
{\vcenter{\offinterlineskip
\halign{\hfil$\scriptscriptstyle##$\hfil\cr>\cr
\noalign{\vskip0.5pt}\simeq\cr}}}}}
\def\grole{\mathrel{\mathchoice {\vcenter{\offinterlineskip\halign{\hfil
$\displaystyle##$\hfil\cr>\cr\noalign{\vskip-1.5pt}<\cr}}}
{\vcenter{\offinterlineskip\halign{\hfil$\textstyle##$\hfil\cr
>\cr\noalign{\vskip-1.5pt}<\cr}}}
{\vcenter{\offinterlineskip\halign{\hfil$\scriptstyle##$\hfil\cr
>\cr\noalign{\vskip-1pt}<\cr}}}
{\vcenter{\offinterlineskip\halign{\hfil$\scriptscriptstyle##$\hfil\cr
>\cr\noalign{\vskip-0.5pt}<\cr}}}}}
\def\leogr{\mathrel{\mathchoice {\vcenter{\offinterlineskip\halign{\hfil
$\displaystyle##$\hfil\cr<\cr\noalign{\vskip-1.5pt}>\cr}}}
{\vcenter{\offinterlineskip\halign{\hfil$\textstyle##$\hfil\cr
<\cr\noalign{\vskip-1.5pt}>\cr}}}
{\vcenter{\offinterlineskip\halign{\hfil$\scriptstyle##$\hfil\cr
<\cr\noalign{\vskip-1pt}>\cr}}}
{\vcenter{\offinterlineskip\halign{\hfil$\scriptscriptstyle##$\hfil\cr
<\cr\noalign{\vskip-0.5pt}>\cr}}}}}
\def\loa{\mathrel{\mathchoice {\vcenter{\offinterlineskip\halign{\hfil
$\displaystyle##$\hfil\cr<\cr\noalign{\vskip1.5pt}\approx\cr}}}
{\vcenter{\offinterlineskip\halign{\hfil$\textstyle##$\hfil\cr<\cr
\noalign{\vskip1.0pt}\approx\cr}}}
{\vcenter{\offinterlineskip\halign{\hfil$\scriptstyle##$\hfil\cr<\cr
\noalign{\vskip0.5pt}\approx\cr}}}
{\vcenter{\offinterlineskip\halign{\hfil$\scriptscriptstyle##$\hfil\cr
<\cr\noalign{\vskip0.5pt}\approx\cr}}}}}
\def\goa{\mathrel{\mathchoice {\vcenter{\offinterlineskip\halign{\hfil
$\displaystyle##$\hfil\cr>\cr\noalign{\vskip1.5pt}\approx\cr}}}
{\vcenter{\offinterlineskip\halign{\hfil$\textstyle##$\hfil\cr>\cr
\noalign{\vskip1.0pt}\approx\cr}}}
{\vcenter{\offinterlineskip\halign{\hfil$\scriptstyle##$\hfil\cr>\cr
\noalign{\vskip0.5pt}\approx\cr}}}
{\vcenter{\offinterlineskip\halign{\hfil$\scriptscriptstyle##$\hfil\cr
>\cr\noalign{\vskip0.5pt}\approx\cr}}}}}
\def\bbbr{\mathbold{R}} 
\def\bbbm{\mathbold{M}}
\def\bbbh{\mathbold{H}}
\def\bbbd{\mathbold{D}}
\def\bbbone{{\mathchoice {\rm 1\mskip-4mu l} {\rm 1\mskip-4mu l}
{\rm 1\mskip-4.5mu l} {\rm 1\mskip-5mu l}}}
\def\bbbc{{\mathchoice {\setbox0=\hbox{$\displaystyle\rm C$}\hbox{\hbox
to0pt{\kern0.4\wd0\vrule height0.9\ht0\hss}\box0}}
{\setbox0=\hbox{$\textstyle\rm C$}\hbox{\hbox
to0pt{\kern0.4\wd0\vrule height0.9\ht0\hss}\box0}}
{\setbox0=\hbox{$\scriptstyle\rm C$}\hbox{\hbox
to0pt{\kern0.4\wd0\vrule height0.9\ht0\hss}\box0}}
{\setbox0=\hbox{$\scriptscriptstyle\rm C$}\hbox{\hbox
to0pt{\kern0.4\wd0\vrule height0.9\ht0\hss}\box0}}}}
\def\bbbe{{\bf E}}
\def\bbbq{{\mathchoice {\setbox0=\hbox{$\displaystyle\rm Q$}\hbox{\raise
0.15\ht0\hbox to0pt{\kern0.4\wd0\vrule height0.8\ht0\hss}\box0}}
{\setbox0=\hbox{$\textstyle\rm Q$}\hbox{\raise
0.15\ht0\hbox to0pt{\kern0.4\wd0\vrule height0.8\ht0\hss}\box0}}
{\setbox0=\hbox{$\scriptstyle\rm Q$}\hbox{\raise
0.15\ht0\hbox to0pt{\kern0.4\wd0\vrule height0.7\ht0\hss}\box0}}
{\setbox0=\hbox{$\scriptscriptstyle\rm Q$}\hbox{\raise
0.15\ht0\hbox to0pt{\kern0.4\wd0\vrule height0.7\ht0\hss}\box0}}}}
\def\bbbt{{\mathchoice {\setbox0=\hbox{$\displaystyle\rm
T$}\hbox{\hbox to0pt{\kern0.3\wd0\vrule height0.9\ht0\hss}\box0}}
{\setbox0=\hbox{$\textstyle\rm T$}\hbox{\hbox
to0pt{\kern0.3\wd0\vrule height0.9\ht0\hss}\box0}}
{\setbox0=\hbox{$\scriptstyle\rm T$}\hbox{\hbox
to0pt{\kern0.3\wd0\vrule height0.9\ht0\hss}\box0}}
{\setbox0=\hbox{$\scriptscriptstyle\rm T$}\hbox{\hbox
to0pt{\kern0.3\wd0\vrule height0.9\ht0\hss}\box0}}}}
\def\bbbs{{\mathchoice
{\setbox0=\hbox{$\displaystyle     \rm S$}\hbox{\raise0.5\ht0\hbox
to0pt{\kern0.35\wd0\vrule height0.45\ht0\hss}\hbox
to0pt{\kern0.55\wd0\vrule height0.5\ht0\hss}\box0}}
{\setbox0=\hbox{$\textstyle        \rm S$}\hbox{\raise0.5\ht0\hbox
to0pt{\kern0.35\wd0\vrule height0.45\ht0\hss}\hbox
to0pt{\kern0.55\wd0\vrule height0.5\ht0\hss}\box0}}
{\setbox0=\hbox{$\scriptstyle      \rm S$}\hbox{\raise0.5\ht0\hbox
to0pt{\kern0.35\wd0\vrule height0.45\ht0\hss}\raise0.05\ht0\hbox
to0pt{\kern0.5\wd0\vrule height0.45\ht0\hss}\box0}}
{\setbox0=\hbox{$\scriptscriptstyle\rm S$}\hbox{\raise0.5\ht0\hbox
to0pt{\kern0.4\wd0\vrule height0.45\ht0\hss}\raise0.05\ht0\hbox
to0pt{\kern0.55\wd0\vrule height0.45\ht0\hss}\box0}}}}

%
% note: changed \sans to \sf for LaTeX
%

\def\bbbz{{\mathchoice {\hbox{$\sf\textstyle Z\kern-0.4em Z$}}
{\hbox{$\sf\textstyle Z\kern-0.4em Z$}}
{\hbox{$\sf\scriptstyle Z\kern-0.3em Z$}}
{\hbox{$\sf\scriptscriptstyle Z\kern-0.2em Z$}}}}

\def\diameter{{\ifmmode\oslash\else$\oslash$\fi}}

\newcommand{\sh}{{\euf s}}
\newcommand{\cbar}{\overline{\bbc}}
\newcommand{\crat}{{\euf c}}
\newcommand{\ccirc}[1]{{\cal C}_{#1}}
\newcommand{\hyp}{\bbbh}
\newcommand{\htwo}{\hyp^2}
\newcommand{\hthree}{\hyp^3}
\newcommand{\etwo}{\bbbe^2}
\newcommand{\ethree}{\bbbe^3}
\newcommand{\reals}{\bbbr}
\newcommand{\comps}{\bbbc}
\def\card{\mathop {\rm card}\nolimits}
\def\grad{\mathop {\rm grad}\nolimits}
%local stuff end here
\newtheorem{Theorem}{Theorem}[section]
\renewcommand{\theTheorem}{\arabic{section}.\arabic{Theorem}}
\newtheorem{Lemma}[Theorem]{Lemma}
\newtheorem{Remark}[Theorem]{Remark}
\newtheorem{Example}[Theorem]{Example}
\newtheorem{Observation}[Theorem]{Observation}
\newtheorem{Proposition}[Theorem]{Proposition}
\newtheorem{Corollary}[Theorem]{Corollary}
\newtheorem{Definition}[Theorem]{Definition}
\newtheorem{Notes}{Notes}
\renewcommand{\theNotes}{}
\newtheorem{Fact}{Fact}
\renewcommand{\theFact}{}
\newtheorem{ParTheorem}{Parametrization Theorem}
\renewcommand{\theParTheorem}{}
\newtheorem{Note}{Note}
\renewcommand{\theNote}{}

%-------------- Author entries --------------------

\title{Simple curves on hyperbolic tori}
%Article title 
\shorttitle{Simple curves on tori} % Shortened version for
					     % headline title 

%% Please enter all acknowledgements here:
\acknowledgements{
The authors would like to thank Peter Sarnak for continued encouragement.
}
% Please uncomment and use appropriate command:
%\author{Igor Rivin}
\twoauthors{Greg McShane}{Igor Rivin}
%\authors{}% Separate each author with a comma and a space.
% Institution:
%% If more than one institution represented, please separate
%% with \\ , i.e.,
%% \institutions{University of Illinois at Chicago, Chicago, IL\\
%% Cornell University, Ithaca, NY}
\institutions{UMPA, Ecole Normale Superieure, 46 Alle\'e d'Italie, 69364 Lyon\\
Melbourne University, Victoria}

{\bf Abstract.} Let $T$ be a once punctured torus, equipped with a complete hyperbolic
metric. Herein, we describe a new approach to
the study of the set ${\cal S}$ of all {\em simple geodesics} on $T.$
We introduce a valuation on the homology  $H_1(T, \bbbz),$ which
associates to each homology class $h$ the
length $\ell(h)$ of the unique simple geodesic homologous to $h,$ and
show that $\ell$ extends to a norm on $H_1(T, \bbbr).$
We analyze the boundary of the unit ball ${\cal B}(\ell)$
and the variation of the area of ${\cal B}(\ell)$ over
the moduli space of $T$.
These results  are applied to obtain sharp asymptotic estimates on the
number of simple geodesics of length less than $L$.

\begin{center}
{\bf Courbes simples dans les tores}
\end{center}
 
{\bf R\'esum\'e.} Soit $T$ un tore trou\'e, muni d'une m\'etrique
hyperbolique compl\`ete, d'aire finie.
Nous pr\'esentons une nouvelle approche  de l'\'etude de l'ensemble
${\cal S}$ de toutes les g\'eod\'esiques ferm\'ees simples (sans points doubles)
de $T$. 
Nous introduisons une application sur l'homologie $H_1(T, \bbbz)$,
qui associe \`a chaque classe  $h\in H_1(T, \bbbz)$ indivisible
la longueur $\ell(h)$ de l'unique g\'eod\'esique simple homologue \`a $h.$
et nous d\'emontrons que $\ell$ s'\'etend en une norme sur $H_1(T,
\bbbr).$ Nous \'etudions la g\'eom\'etrie de la sph\`ere $\partial {\cal
B}(\ell)$ et la variation de l'aire de ${\cal B}(\ell)(T)$ sur l'espace
des modules. On utilise ces r\'esultats pour donner des estimations
asymptotiques  du nombre de g\'eod\'esiques ferm\'ees simples de longueur
inf\'erieure \`a $L$. 
 
{\bf Version francaise abr\'eg\'ee}
Soit $T$ un tore trou\'e, muni d'une m\'{e}trique hyperbolique compl\`{e}te
d'aire finie. Dans cet article nous d\'{e}crivons une m\'{e}thode pour
\'etudier l'ensemble ${\cal S}$ des {\em lacets g\'{e}od\'{e}siques simples}
sur $T.$ L'ensemble ${\cal S}$ est un objet assez int\'{e}ressant
(v. \cite{haas} pour des applications arithm\'{e}tiques). 

Soit $S$ une surface hyperbolique ferm\'{e}e, ou avec un seul cusp. Une
{\em multi-courbe} $m$ sur $S$ est une application continue d'une vari\'{e}t\'{e} de
dimension $1$ dans $S.$ La {\em longueur} d'une multi-courbe est d\'{e}finie
comme la somme des longueurs des images de composantes connexes de $M.$ Nous
dirons que $m$ est plong\'{e}e lorsque l'image de $m$ est une r\'{e}union
disjointe de lacets simples -- l'application $m$ n'est pas forc\'ement
injective. Notre premier r\'{e}sultat est que pour chaque classe $h \in
H_1(S, \bbbz),$ il y a une multi-courbe $m$ de longueur minimale,
homologue \`a $h.$ Cette multi-courbe $m$ est plong\'{e}e, et toutes
les composantes connexes de $m$ sont des lacets g\'{e}od\'{e}siques.  

Dans le cas sp\'ecial o\`{u} $S$ est un tore trou\'{e}, nous savons qu'il
y a un seul lacet g\'{e}od\'{e}sique simple dans chaque $h \in H_1(S,
\bbbz),$ quand $h$ est indivisible, et que deux lacets
g\'{e}od\'{e}siques simples quelconques se rencontrent. Donc, il y a une seule
multi-courbe minimale $m_h$ dans chaque classe d'homologie $h$ -- si $h$
est indivisible, $m_h$ est l'unique lacet g\'{e}od\'{e}sique simple homologue
\`{a} $h$ et, autrement, $m_h$ rev\^et un lacet g\'{e}od\'{e}sique simple.

Maintenant, nous pouvons d\'{e}finir une valuation $\ell: H_1(T, \bbbz)
\rightarrow \bbbr,$ telle que $\ell(h)$ est la longueur de la
(seule) multi-courbe homologue \`{a} $h.$ Nous pouvons d\'{e}montrer
(en utilisant la discussion pr\'{e}c\'{e}dente) que $\ell$ s'\'{e}tend
en  une  {\em norme} sur $H_1(T, \bbbr) \simeq \bbbr^2.$ Soit ${\cal
B}_\ell$ la boule unit\'e de cette norme -- ${\cal B}_\ell$ est un
convexe sym\'{e}trique dans le plan $\bbbr^2.$ La sph\`{e}re unit\'e
$\partial {\cal B}_\ell$ est une repr\'{e}sentation de l'espace projectif des
laminations mesur\'{e}es du tore trou\'{e}. 

La g\'{e}om\'{e}trie de la boule ${\cal B}_\ell$ est assez int\'{e}ressante. Le
th\'eor\`eme suivant sugg\`ere que la sph\`{e}re $\partial {\cal B}_\ell$ est 
vraiment irr\'{e}guli\`{e}re:

{\bf Th\'eor\`eme.}
{
\em
Le bord de la boule unit\'e, $\partial {\cal B}_\ell,$ 
a un coin en tous les points de pente irrationnelle.
L'angle ext\'erieur au coin de pente $p/q$
d\'ecro\^{\i}t exponentiellement comme fonction de $\max(p,q)$.
En tout point de pente irrationnelle le bord est infiniment plat.
}

Le terme ``infiniment plat'' est d\'efini de la
mani\`ere suivante: 

Apr\`es une isom\'etrie du plan euclidien, nous pouvons
supposer que le point $p$ est l'origine du plan, une droite de contact  
\`a ${\cal B}_\ell$ dans $p$ est la droite $y=0$ et la boule ${\cal
B}_\ell$ est situ\'ee dans le demi-espace $y\geq 0.$ Il existe un
voisinage $|x| < \delta_0,$ tel que le bord de ${\cal B}_\ell$ est un
graphe d'une fonction convexe $f.$ Nous disons que $\partial {\cal
B}_\ell$ est {\em infiniment plate} en $p$, si pour tout entier $N,$
il existe un $0<\delta<\delta_0,$ tel que $f(x) \leq |x|^N,$ quand
$|x| < \delta.$ 

Nous appliquons ces r\'{e}sultats pour calculer le nombre asymptotique
de lacets g\'{e}od\'{e}siques de longueur born\'{e}e par $L$ sur le
tore trou\'{e} $T.$ Nous d\'emontrons que 
$$N_L = \frac{L^2}{\Area {\cal B}_\ell} + O(L \log L).$$

Ce r\'esultat est assez simple modulo la discussion pr\'ec\'edente -- la
quantit\'e $N_L$ est \'egale au nombre des points entiers {\em
indivisibles} -- $(m, n)$ est indivisible quand le rayon de l'origine \`a
$(m, n),$ ne contient aucun point entier plus petit que $(m, n);$
Hardy et Wright utilisent le terme ``visible'' -- 
dans $L {\cal B}_\ell.$ Sans l'hypoth\`ese d'indivisibilit\'e, l'erreur
est $O(L);$ l'estimation pour les points s'ensuit apr\`es une inversion
de M\"obius. C'est un fait remarquable que cette estimation est
presque exacte:

{\bf Th\'eor\`eme.} {\em Il y a une partie dense de l'espace des modules pour
laquelle l'erreur ci-dessus est plus grande que $O(L).$ }

Nous \'{e}tudions aussi la g\'{e}om\'{e}trie de
${\cal B}_\ell(T)$ comme une fonction de $T$ dans l'espace des modules
${\cal M}_{1, 1}.$ Nous d\'{e}montrons que l'aire de ${\cal B}_\ell(T)$
tend vers l'infini quand $T$ s'approche des cusps de l'espace de
modules. 

La structure hyperbolique la plus sym\'etrique, et d'int\'er\^et
arithm\'etique maximal est {\em le tore trou\'e modulaire} -- 
le tore trou\'e $T_m$
unique tel que $T_m = \Gamma \backslash {\bf H}^2,$ o\`u $\Gamma$ est un
sous-groupe de $SL(2, \bbbz).$ Nous croyions que:

{\bf Conjecture 1.} L'aire de ${\cal B}_\ell(T)$ est minimale pour $T=
T_m$ (dans cet cas, l'aire est \'egale \`a $A_m = 1.08...$.

et aussi

{\bf Conjecture 2.} Pour $T_m,$ l'estimation n'est pas g\'en\'erique.
C'est-\`a-dire:
$$
N_L(T_m) = L^2/A_m + O(L^{\frac12 + \epsilon}),
$$
pour tout $\epsilon,$ comme $L \rightarrow \infty$.

Nous croyons que la deuxi\`eme conjecture est tr\`es difficile.

\intro

Let $T$ be a torus with one puncture, equipped with a complete hyperbolic
metric of finite area. In this paper we describe a new approach to
the study of the set ${\cal S}$ of all {\em simple geodesics} (that is,
those without  self-intersections) on $T.$ The set ${\cal S}$ is of
considerable interest in both geometry and number theory (see \cite{haas}). 

The plan of the paper is as follows. 
In sections \ref{minisec} and \ref{norm}, we introduce a valuation on
the homology group $H_1(T, \bbbz),$ which associates to each homology
class $h$ the 
length $\ell(h)$ of the unique simple geodesic homologous to $h,$ and
show that $\ell$ extends to a norm on $H_1(T, \bbbr).$ In
section \ref{norm} we investigate the geometry of the norm $\ell,$
by analyzing the boundary of the unit ball${\cal B}(\ell).$ In section
\ref{modvar} we analyze the behavior of $\ell$ as $T$ varies over the
moduli space ${\cal M}^{1, 1}$ of finite area complete hyperbolic
metrics on the punctured torus. Finally, in section \ref{asymp} we
apply these results to obtain asymptotic estimates on the number of
simple geodesics with length bounded above by $L.$

\section{Minimal multicurves}
\label{minisec}

Let $S$ be a hyperbolic surface of finite volume with at most one
cusp. We define a {\em multicurve} $m$ on $S$ to be a map from a (not
necessarily connected) 
$1$-manifold $M$ to $S.$ We define the {\em length} of $m$ to be the
sum of the lengths of the images of components of $M.$ We say that a
multicurve $m$ is {\em embedded} if the image of $m$ is the union of
simple closed curves $\gamma_1, \dots, \gamma_k$ on $S.$ Note that the
map $m$ may cover some of the components $\gamma_i$ multiple times, so
this does not coincide with the usual meaning of embedding. A
multicurve defines a singular chain, which, in turn, defines a
homology class in $H_1(S, \bbbz).$

\thm{Theorem}{\label{minimal}
Let $h \in H_1(S, \bbbz)$ be a non-trivial homology class. There
exists a multicurve $m$ representing $h$ of minimal length, and $m$ is
embedded, with all components geodesic. 
}

\thm{Corollary}{\label{mintor}
Let $T$ be a punctured torus equipped with a hyperbolic structure.
Then, the shortest multicurve representing a non-trivial homology
class $h$ is a simple closed geodesic if $h$ is a primitive homology
class (that is, not a multiple of another class), and a multiply
covered geodesic otherwise. In addition, the shortest multicurve
representing $h$ is unique.}

The reader can find the proofs of these results in \cite{mcriv}.

\section{A norm on homology}
\label{norm}

For each primitive homology class $h \in H_1(T, \bbbz),$ define
$\ell(h)$ to be the length of the unique simple geodesic $\gamma_h
\sim h.$ The uniqueness is guaranteed by Corollary \ref{mintor} of the
previous section. For a non-primitive homology class $(m, n),$ define $\ell(h)
= \gcd (m, n) \ell(h/(\gcd(m, n))).$ In this section we will show that
$\ell$ can be extended to a {\em norm} on $H_1(T, \bbbr) \simeq
\bbbr^2.$ 

The results of section \ref{minisec} imply, in particular, that $\ell$
satisfies a triangle inequality on $H_1(T, \bbbz),$ 

\demo{Note} The triangle inequality is {\em strict} (that is, if
$\ell(n, m) + \ell(p, q) = \ell(n+p, m+q),$ then $(n, m)$ and $(p, q)$
are both multiples of the same class). 
\enddemo

We can now extend $\ell$ to $H_1(T, \bbbr);$ first to $H_1(T, \bbbq)$
by linearity, and then to $H_1(T, \bbbr)$ by continuity (which is an
immediate consequence of the triangle inequality). It is, furthermore,
clear that $\ell$ still satisfies a triangle inequality, and is, thus,
a {\em semi-norm}. To show that it is a norm, requires some study of
simple geodesics from a group-theoretic standpoint, and some
elementary techniques of geometric group actions.
The reader is referred to \cite{mcriv} for the details.

In the figure, the reader can find a computer-generated picture of the
unit ball ${\cal B}$ of the norm $\ell$ for $T$ -- the modular torus (the only one
corresponding to a subgroup of $SL(2, \bbbz)$, whose fundamental group
is generated by the matrices $\pmatrix{1 & 1 \cr 1 & 2}$ and
$\pmatrix{1 & -1 \cr -1 & 2}$). 

\vspace{3in}
\epsfysize=4in
\epsffile{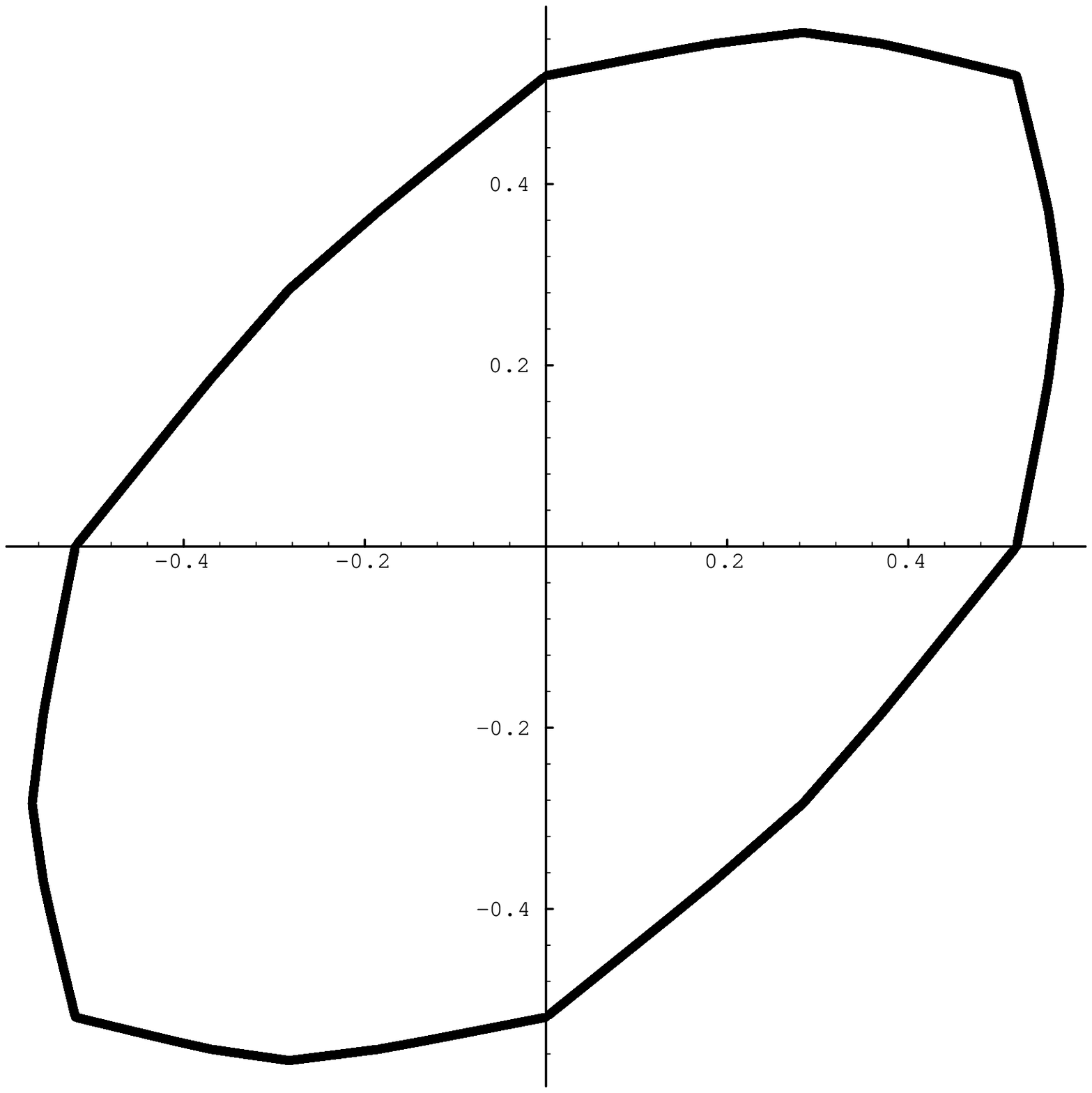}

%figure here
Several very accurate computer generated renderings of the unit ball
were made by the authors. To the casual observer first seeing
these pictures  ${\cal B}$ appears polygonal, 
but it is not hard to see that this is far from the truth -- since
the triangle inequality is strict, no three points with rational slope
are collinear. To describe the true geometry of ${\cal B}$ (for any
complete hyperbolic structure of finite area), we first need to define
some terms.

Let $B$ be a convex body in the plane, and let $p \in \partial B.$
After applying a euclidian isometry we can 
suppose that the point $p$ is at the origin in the plane and that
$B$ lies in the half space  $y \geq 0$. There
is a neighborhood $|x| < \delta_0$, such that the the boundary
of $\partial B$ is a graph of a convex function $f$. 
We say that $\partial B$ is
{\em flat to order $n$} at $p$ if there exists a $\delta,$
$0<\delta<\delta_0,$ such that $f(x) \leq |x|^n$, when $|x| < \delta.$
We say that $\partial B$ is {\em flat to infinite order} at $p$, if it
is flat to order $n,$ for every $n>0.$

Now we can state:
\thm{Theorem}
{
\label{ballbdry}
The boundary of the unit ball ${\cal B}_\ell$ has a corner at each
point of rational slope. The exterior angle at such a corner at the
slope $p/q$ decreases exponentially as a function of $\max(p, q)$ (the {\em
height}). At a point of irrational slope $\theta,$ the boundary of
${\cal B}_\ell$ is flat, to infinite order.
}

The proof of Theorem \ref{ballbdry} is somewhat involved (and
presented in full in \cite{mcriv2}). Among the main ingredients are the
following: 

{\bf Fricke trace relation.} This says that for every pair of elements
$A,$ $B$ of $SL(2, \bbbr)$ 
$$
2 + \tr[A, B] = \tr^2 A + \tr^2 B + \tr^2 AB - \tr A \tr B \tr AB,
$$
where $[A, B] (= ABA^{-1}B^{-1})$ is the commutator of $A$ and $B.$ (See, eg, \cite{fenchel}).

When $A$ and $B$ are associated generators of the fundamental group of
a cusped torus, their commutator is a parabolic element (a simple loop
around the cusp), and so the left hand side of this relation vanishes.

{\bf Elliptic involution.} The punctured torus is elliptic, and has
four {\em Weierstrass points}, preserved by the elliptic
involution. One of these is the cusp, the other three have the
property that every simple closed geodesic passes through two of them, and
further if $\alpha$ passes through $W_1, W_2$ and $\beta$ passes
through $W_1, W_3,$ then $\alpha \beta$ passes through $W_2$ and
$W_3.$ 

The trace relation can be used to compute the angle between the geodesics
$\alpha$ and $\beta$ at the point $W_1$ in terms of their
lengths, as follows: 

\thm{Theorem}
{
\label{sinsin}
$$\sinh\ell(\alpha) \sinh \ell(\beta) \sin\angle \alpha, \beta = 1.$$
}

This identity can then be used to show that the angle between $\alpha$
and $\beta$ goes to $0$ (or $\pi$) exponentially fast with the greater
of the two lengths, and this, in turn leads to a proof that the
triangle inequality is exponentially close to equality, again, as a
function of the greater of the two lengths. This, in turn, is the main
tool in the proof of Theorem \ref{ballbdry}.

\section{Variation over moduli space}
\label{modvar}

In this section we describe the variation of ${\cal B}_\ell$ over the
moduli space of complete finite area hyperbolic metrics on the
punctured torus. First, we note that points on
$\partial {\cal B}_\ell$ correspond to {\em projective measured
laminations} on the punctured torus -- thus ${\cal B}_\ell$ can be
viewed as an explicit picture of the projective lamination space
${\cal P}{\cal L}.$ It should be noted that a somewhat similar
combinatorial picture of that space is constructed in
\cite{thurmin}. It is not hard to show that the topology of
convergence in $\ell$ is exactly the topology on ${\cal P}{\cal L}.$
Thus we can use Kerckhoff's results \cite{earth} to show the following:
\thm{Theorem}
{
\label{analvar}
Every point of $\partial {\cal B}_\ell(T)$ varies analytically over moduli space.
}
\thm{Corollary}
{
\label{analarea}
The area of ${\cal B}_\ell(T)$ varies analytically over moduli space.
}
\thm{Theorem}
{
\label{cuspvar}
The area of ${\cal B}_\ell(T)$ goes to infinity as $T$ approaches the
cusps of the moduli space.
}

\demo{Proof} The cusps of the moduli spaces are characterized by the
condition that the length of some simple geodesic $\gamma$ approaches
to $0.$ Pick a point $T_0$ of moduli space, pick a simple geodesic
$\gamma,$ and let $\gamma^\prime$ be the {\em shortest} associated
generator of the fundamental group. 
Since $\gamma^\prime$ is shortest, it follows that  $\ell(\gamma
\gamma^\prime) \geq \ell(\gamma^\prime)$ and $\ell(\gamma^\prime
\gamma^{-1}) \geq \ell(\gamma^\prime).$ Keeping in mind that $\cosh
\ell(\gamma)/2 = \tr \gamma/2,$ (where we abuse 
notation by using $\gamma$ both for the geodesic and the corresponding
element of $SL(2, \bbbr),$ the above inequality and the Fricke trace
relation imply:
\begin{equation}
\label{eq1}
\cosh \frac{\ell(\gamma^\prime)}2 \leq \cosh \frac{\ell (\gamma^\prime
\gamma^{-1})}2 = \cosh \frac{\ell(\gamma^\prime)}2 \cosh \frac{\ell(\gamma)}2
\left(1 - \sqrt{1-\frac{1}{\cosh^2 \frac{\ell(\gamma^\prime)}2} - 
\frac{1}{\cosh^2 \frac{\ell(\gamma)}2}}\right).
\end{equation}
A computation transforms the above into
\begin{equation}
1-\frac{1}{\cosh^2 \frac{\ell(\gamma)}{2}} \geq \frac{1}{\cosh^2
\frac{\ell(\gamma^\prime)}2} \geq  \frac{2}{\cosh
\frac{\ell(\gamma)}2} - \frac{2}{\cosh^2 \frac{\ell(\gamma)}2}, 
\end{equation}
where the first inequality stems from the positivity of the expression
inside the radical in equation \ref{eq1}

When $\gamma$ is very short, we can expand $\cosh
\frac{\ell(\gamma)}2$ in a Taylor series, to obtain
$$\left| \frac{1}{\cosh^2 \frac{\ell(\gamma^\prime)}2} -
\left(\frac{\ell(\gamma)}2\right)^2 \right| = O(\ell^4(\gamma)) \Longrightarrow
 \ell(\gamma^\prime) \sim |\log \ell(\gamma)|.$$
Since 
$$\Area {\cal B}_\ell \geq \frac{1}{\ell(\gamma)
\ell(\gamma^\prime)},$$ it follows that 
$\lim_{\ell(\gamma)\rightarrow 0} \Area {\cal B}_\ell =\infty.$
\enddemo

The above proof actually shows a little more -- by combining the
estimates with Theorem \ref{sinsin}, we have
\thm{Corollary}
{
The angle between a very short simple geodesic $\gamma$ and the shortest
conjugate geodesic $\gamma^\prime$ approaches $\pi/2,$ as $\ell
(\gamma) \rightarrow 0.$
}

\demo{Definition} The {\em systole} $\sys T$ is the shortest closed
geodesic on $T.$
\enddemo
A modification of the proof of Theorem \ref{cuspvar} can be used to
obtain the following result:

\thm{Theorem}
{
\label{syst}
The modular torus is {\em isosystolic}:
$$\length \sys T_{\bmod} = \max_{T \in {\cal M}^{1, 1}} \length \sys T.$$
}

We would like to propose the following conjecture:

\thm{Conjecture}
{
The area of ${\cal B}_\ell(T)$ is minimized when $T$ is the modular
torus.
}

\section{Number of simple geodesics with length bounded above}
\label{asymp}

The results of section \ref{norm} tell us that the number of simple
geodesics of length bounded by $L$ is just the number of primitive
lattice points in $L {\cal B}_\ell,$ which, by the usual trivial
lattice point estimate combined with M\"obius inversion tells us that

\thm{Theorem}
{
\label{latest}
The number of simple geodesics of length bounded by $L$ is 
$$N(L)=\frac{1}{\zeta(2)}\Area({\cal B}_\ell) L^2 + O(L \log L),$$
where $\zeta$ is the Riemann zeta function.
}

A similar estimate for the modular torus was obtained (in
a totally different language, and by number-theoretic methods) by
D.~Zagier in \cite{zagier}. He conjectured that the error bound was
essentially sharp. Indeed, we can show:

\thm{Theorem}
{
\label{esharp}
The error term in Theorem \ref{latest} is {\em at least} O(L) for a
dense set in moduli space.
}

The difference between $O(L)$ and $O(L \log L)$ is caused by the
difference between counting all lattice points and primitive lattice
points. To show theorem \ref{esharp}, it is enough to observe that
since the direction of the support line at any fixed slope varies
analytically over moduli space, its slope will be rational at a point
of irrational slope at a dense set of points in moduli space. Call one
of these points $T_r.$ Now, consider a convex set $S$ in the plane,
such that there is a $p \in \partial S,$ such that the curvature of
$\partial S$ at $p$ vanishes to infinite order, and the slope of the
normal is rational.
It can be shown (see \cite{guillemin}) that
the number of lattice points in $t S$ is $ t^2 \Area S + O(t),$
and the error term is {\em sharp}.
It follows from the proof of that fact, that the error term for the
{\em primitive} lattice point problem is at least linear. Thus, the
error term for $T_r$ in theorem \ref{esharp} is seen to be somewhere
between $O(L)$ and 
$O(L \log L).$ While it is not at all clear whether the error term is
better than linear anywhere in moduli space, we put forward the
following conjecture:

\thm{Conjecture}
{
The error term in Theorem \ref{esharp} for the {\em modular} torus is
of order $O(L^{\frac12+\epsilon}).$
}

This is contrary to Zagier's conjecture which was also stated for the
modular torus.

\end{document}